\documentclass[11 pt, psamsfonts]{amsart}
\usepackage{epic}
\usepackage{eepicemu}

\usepackage{amssymb}

\def\ignore#1{\relax}

\def\twlrm{}
\newdimen\ybox\ybox=5pt
\setbox249=\vbox{\hbox{\vrule{\vbox to\ybox{}}\kern\ybox}\hrule}
\setbox248=\hbox{\copy249}

\newcount\n
\def\R#1,{%
\n=#1%
\nointerlineskip
\hbox{%
\loop\ifnum\n>0\copy248\advance\n by-1\repeat%
\vrule}%
}

\def\YD#1.{\,\vcenter{\hrule#1}\,}

\def\nat{{\mathbb N}}

\def\C{{\mathbb C}}

\def\la{\lambda}

\def\one{\mathbf 1}

\def\v{\vskip 2.5mm}

\def\inv{^{-1}}
\def\ignore#1{\relax}
\def\one{{\bf 1}}

\def\TL{{\rm{TL}}}
\def\Tr{{\rm{Tr}}}
\def\Hom{{\rm{Hom}}}
\def\End{{\rm{End}}}
\def\id{{\rm{id}}}
\def\tr{{\rm{tr}}}
\def\Ne{{\rm Neg}}
\def\Neg{{\rm Neg}}
\def\ve{\varepsilon}
\def\glm{g_{\la, \mu}}
\def\al{\alpha}

{\theoremstyle{plain}
\newtheorem{theorem}{Theorem}[section]
}

{\theoremstyle{plain}
\newtheorem{proposition}[theorem]{Proposition}
}

{\theoremstyle{plain}
\newtheorem{corollary}[theorem]{Corollary}
}

{\theoremstyle{plain}
\newtheorem{lemma}[theorem]{Lemma}
}

{\theoremstyle{plain}

}

{\theoremstyle{definition}

}

{\theoremstyle{definition}

}

{\theoremstyle{remark}

}

{\theoremstyle{remark}

}

\numberwithin{equation}{section}

\renewcommand{\labelenumi}{{ \theenumi.}}

\begin{document}

\title[Temperley-Lieb Ideals]
{Ideals in the Temperley-Lieb Category}

\author{Frederick M. Goodman and Hans Wenzl}

\address{Department of Mathemematics\\ University of Iowa\\ Iowa City,
Iowa}

\email{goodman@math.uiowa.edu}

\address{Department of Mathematics\\ University of California\\ San
Diego,
California}

\email{wenzl@brauer.ucsd.edu}

\thanks{We are grateful to Michael Freedman for bringing the question of tensor
ideals in the Temperley-Lieb category  to our attention and for allowing
us to present the proof as an appendix to his paper.}

\date{written April, 2001}

\bigskip

\maketitle

This note will appear as an appendix to the paper of 
 Michael Freedman,
{\em A magnetic model with a possible Chern-Simons
phase}~\cite{Freedman}; it may, however, be read independently of
~\cite{Freedman}.  The purpose of this paper  is to prove
the following result.

\begin{theorem}
For a generic value of the parameter, the Temperley-Lieb category has no
non-zero, proper tensor ideal.  When the parameter $d$ is equal to 
$2\cos(\pi/n)$ for some $n \ge 3$, then the Temperley-Lieb category has
exactly one non-zero, proper ideal, namely the ideal of negligible
morphisms.
\end{theorem}

Our notation in the appendix differs slightly from that in
~\cite{Freedman}. We write $t$ instead of $-A^2$,  $T_n$ for the
Temperley-Lieb algebra with $n$ strands, and $TL$ for the Temperley-Lieb
category. We trust that this notational variance will not cause the
reader any difficulty.

\section{The Temperley-Lieb Category}

\subsection{The Generic Temperley Lieb Category}

Let $t$ be an indeterminant over $\C$, and let $d = (t + t\inv)$.    The
{\em
generic Temperley Lieb category}  TL is a strict tensor categor whose
{\em
objects}  are elements of
$\nat_0 = 
\{0, 1, 2, \dots\}$.
The set of {\em morphisms} $\Hom(m,n)$ from $m$ to $n$ is a $\C(t)$
vector
space described as follows:

\smallskip
If $n-m$ is odd, then $\Hom(m,n)$ is the zero vector space.

For $n-m$ even, we first define
$(m,n)$--TL diagrams,  consisting of:
\smallskip\noindent
\begin{enumerate}
\item  a closed rectangle $R$  in the plane with two opposite edges
designated as top and bottom.
\item $m$  marked points (vertices) on  the top edge and $n$ marked
points on the bottom edges.
\item $(n+m)/2$ smooth curves (or ``strands") in $R$ such that for each
curve
$\gamma$,
$\partial
\gamma = \gamma \cap \partial R$ consists of two of the $n+m$  marked
points, and such that the curves are pairwise non-intersecting.
\end{enumerate}

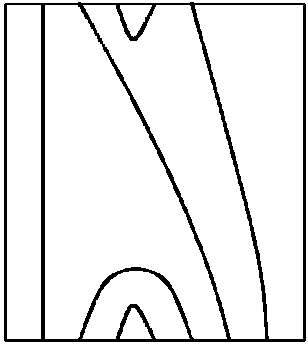
\begin{figure}[ht]
\setlength{\unitlength}{0.00087489in}
\begingroup\makeatletter\ifx\SetFigFont\undefined%
\gdef\SetFigFont#1#2#3#4#5{%
  \reset@font\fontsize{#1}{#2pt}%
  \fontfamily{#3}\fontseries{#4}\fontshape{#5}%
  \selectfont}%
\fi\endgroup%
{\renewcommand{\dashlinestretch}{30}
\begin{picture}(1824,2084)(0,-10)
\path(12,2047)(1812,2047)(1812,22)
    (12,22)(12,2047)
\thicklines
\path(462,2047)(463,2046)(464,2044)
    (466,2040)(470,2034)(475,2025)
    (482,2013)(490,1998)(500,1981)
    (512,1960)(526,1936)(541,1909)
    (557,1880)(575,1849)(594,1816)
    (613,1781)(634,1745)(654,1707)
    (676,1668)(698,1628)(720,1587)
    (742,1545)(765,1502)(788,1457)
    (812,1411)(836,1364)(861,1315)
    (886,1264)(912,1211)(939,1156)
    (966,1099)(993,1041)(1020,982)
    (1047,922)(1076,855)(1104,791)
    (1129,730)(1153,672)(1174,618)
    (1194,567)(1212,520)(1228,476)
    (1242,434)(1256,395)(1268,358)
    (1280,322)(1290,289)(1300,256)
    (1309,226)(1317,197)(1324,170)
    (1331,145)(1337,121)(1343,101)
    (1347,82)(1351,66)(1355,53)
    (1357,42)(1359,34)(1360,28)
    (1361,25)(1362,23)(1362,22)
\path(237,2047)(237,22)
\path(462,22)(463,25)(465,32)
    (470,45)(476,63)(484,86)
    (494,114)(505,144)(517,175)
    (529,205)(541,234)(552,261)
    (564,286)(575,308)(586,328)
    (596,345)(607,360)(618,374)
    (630,386)(642,397)(655,407)
    (668,416)(682,424)(697,431)
    (713,437)(729,442)(746,446)
    (764,449)(781,451)(800,451)
    (818,451)(835,449)(853,446)
    (870,442)(886,437)(902,431)
    (917,424)(931,416)(944,407)
    (957,397)(969,386)(981,374)
    (992,360)(1003,345)(1013,328)
    (1024,308)(1035,286)(1047,261)
    (1058,234)(1070,205)(1082,175)
    (1094,144)(1105,114)(1115,86)
    (1123,63)(1129,45)(1134,32)
    (1136,25)(1137,22)
\path(1137,2047)(1138,2045)(1139,2041)
    (1141,2033)(1144,2021)(1149,2004)
    (1155,1981)(1163,1952)(1173,1917)
    (1184,1877)(1196,1832)(1210,1783)
    (1224,1729)(1240,1673)(1255,1614)
    (1272,1555)(1288,1495)(1304,1435)
    (1320,1376)(1336,1318)(1351,1262)
    (1365,1208)(1379,1157)(1392,1107)
    (1405,1060)(1416,1016)(1427,974)
    (1438,934)(1448,896)(1457,860)
    (1465,826)(1473,794)(1481,764)
    (1488,734)(1494,707)(1501,680)
    (1506,654)(1512,629)(1520,592)
    (1528,556)(1535,521)(1541,488)
    (1546,456)(1552,425)(1556,393)
    (1561,362)(1564,330)(1568,298)
    (1571,266)(1574,233)(1577,201)
    (1579,170)(1581,139)(1583,111)
    (1584,87)(1585,66)(1586,49)
    (1586,36)(1587,28)(1587,24)(1587,22)
\path(687,2047)(687,2046)(689,2041)
    (693,2029)(699,2009)(708,1984)
    (717,1956)(728,1928)(737,1903)
    (746,1881)(755,1865)(763,1852)
    (770,1843)(777,1839)(785,1837)
    (792,1839)(800,1843)(809,1852)
    (818,1865)(829,1881)(841,1903)
    (854,1928)(868,1956)(882,1984)
    (894,2009)(904,2029)(909,2041)
    (912,2046)(912,2047)
\path(687,22)(687,23)(689,28)
    (693,40)(699,60)(708,85)
    (717,113)(728,141)(737,166)
    (746,188)(755,204)(763,217)
    (770,226)(777,230)(785,232)
    (792,230)(800,226)(809,217)
    (818,204)(829,188)(841,166)
    (854,141)(868,113)(882,85)
    (894,60)(904,40)(909,28)
    (912,23)(912,22)
\end{picture}
}
\caption{A (5,7)--Temperley Lieb Diagram}
\label{TL diagram}
\end{figure}

Two such diagrams are {\em equivalent} if they induce the same pairing
of
the
$n+m$ marked points.
 $\Hom(m,n)$ is defined to be the $\C(t)$  vector space
with basis the set of equivalence classes of $(m,n)$--TL diagrams; we
will refer to equivalence classes of diagrams simply as diagrams.

The composition of morphisms is defined first on the level of diagrams.
The
composition
$ba$ of an
$(m,n)$--diagram
$b$ and an
$(\ell,m)$--diagram
$a$  is defined by the following steps:
\smallskip\noindent

\begin{enumerate}
\item
 Juxtapose the rectangles of $a$ and $b$, identifying the
bottom edge of $a$ (with its $m$  marked points) with the top edge of
$b$
(with its $m$ marked points).
\item Remove from the resulting rectangle any closed loops in its
interior.  The result is a $(n,\ell)$--diagram $c$.
\item The product $ba$ is $d^r c$, where $r$ is the number of
closed loops removed.
\end{enumerate}

The composition product evidently respects equivalence of diagrams, and
extends
uniquely to a bilinear product
$$\Hom(m,n)
\times
\Hom(\ell,m)
\longrightarrow  \Hom(\ell,n),$$ hence to a linear map
$$\Hom(m,n) \otimes
\Hom(\ell,m)
\longrightarrow  \Hom(\ell,n).$$

The {\em tensor product of objects} in $\TL$ is given by $n\otimes n' =
n+n'$. The {\em tensor product of morphisms}  is defined by horizontal
juxtposition.  More exactly,  The tensor $a \otimes b$ of
an
$(n,m)$--TL diagram
$a$ and an $(n',m')$--diagram $b$ is defined by horizontal juxtposition
of
the diagrams, the result being an $(n+n',m+m')$--TL diagram.

The tensor product extends uniquely to a bilinear product $$\Hom(m,n)
\times
\Hom(m',n')
\longrightarrow  \Hom(m+m',n+n'),$$ hence to a linear map
$$\Hom(m,n) \otimes
\Hom(m',n')
\longrightarrow  \Hom(m+m',n+n').$$

For each $n \in \nat_0$, $T_n := \End(n)$ is a $\C(t)$--algebra,
with the composition product.  The identity $1_n$ of $T(n)$ is the
diagram with $n$ vertical (non-crossing) strands.  We have canonical
embeddings of $T_n$ into $T_{n+m}$ given by $x \mapsto x \otimes 1_m$.
If
$m>n$ with $m-n$ even, there also exist obvious embeddings of
$\Hom(n,m)$ and $\Hom(m,n)$ into $T_m$ as follows: If $\cap$ and $\cup$
denote the morphisms in $\Hom(0,2)$ and $\Hom(2,0)$, then
we have linear embeddings
$$a\in\Hom(n,m)\mapsto a\otimes \cup^{\otimes (m-n)/2}\in T_m$$
and 
$$b\in\Hom(m,n)\mapsto b\otimes \cap^{\otimes (m-n)/2}\in T_m.$$

Note that these maps have left inverses which are given by
premultiplication
by an
element of $\Hom(n, m)$ in the first case, and postmultiplication by an
element of
$\Hom(m,n)$ in the second.  Namely,
$$
a = d^{-(m-n)/2} (a\otimes \cup^{\otimes (m-n)/2}) \circ
(\one_n \otimes \cap^{\otimes
(m-n)/2})  
$$
and
$$
b = d^{-(m-n)/2} 
(\one_n \otimes \cup^{\otimes
(m-n)/2})   \circ (b\otimes \cap^{\otimes (m-n)/2})
$$

By an {\em ideal} $J$ in TL  we shall mean a vector
subspace of  $\bigoplus_{n,m} \Hom(n,m)$ which is closed under
composition
and tensor product with arbitrary morphisms.  That is, if $a, b$ are
composible morphisms, and one of them is in $J$, then the composition
$ab$ is in $J$; and if $a, b$ are any morphisms, and one of them is in
$J$,
then the tensor product $a \otimes b$ is in $J$.

Note that any ideal is closed under the embeddings described just above,
and
under
their left inverses.

\subsection{Specializations and evaluable morphisms.}  For any $\tau \in
\C$,
we define the specialization $\TL(\tau)$ of the Temperley Lieb category
at
$\tau$, which is obtained by replacing the indeterminant $t$ by $\tau$.
More
exactly,  the objects of $\TL(\tau)$ are again elements of $\nat_0$, the
set
of morphisms
$\Hom(m,n)(\tau)$ is the $\C$--vector space with basis the set of
$(m,n)$--TL diagrams, and the composition rule is as before, except that
$d$ is replaced by $d(\tau) = (\tau + \tau\inv)$.  Tensor products
are defined as before.
$T_n(\tau) := \End(n)$ is
a complex algebra, and $x \mapsto x \otimes 1_m$ defines a canonical
embedding of $T_n(\tau)$ into $T_{n+m}(\tau)$.  One also has embeddings
$\Hom(m,n) \rightarrow T_n$ and $\Hom(n,m) \rightarrow T_n$, when $m <
n$, as before.   An ideal in $\TL(\tau)$ again means a subspace of
$\bigoplus_{n,m} \Hom(n,m)$ which is closed under composition
and tensor product with arbitrary morphisms.

Let $\C(t)_\tau$ be the ring of rational functions without pole at
$\tau$.
The set of {\em evaluable} morphisms in $\Hom(m,n)$ is the
$\C(t)_\tau$--span
of the basis of $(n,m)$--TL diagrams.  Note that the composition and
tensor product of evaluable morphisms are evaluable.  We have an {\em
evaluation map} from the set of evaluable morphisms to morphisms of
$\TL(\tau)$ defined by
$$a=\sum s_j(t)a_j \mapsto a(\tau)=\sum s_j(\tau)a_j,$$
where the $s_j$ are in $\C(t)_\tau$, and the $a_j$ are TL-diagrams.
We write $x \mapsto x(\tau)$ for the evaluation map.
The evaluation map is a homomorphism for the composition and tensor
products.
In particular, one has a $\C$--algebra homomorphism from the algebra
$T_n^\tau$ of evaluable endomorphisms of $n$ to the algebra $T_n(\tau)$
of
endomorphisms of $n$ in $\TL(\tau)$.

The principle of constancy of dimension  is an
important tool for analyzing the specialized categories $\TL(\tau)$.  We
state it
in the form which we need here:

\begin{proposition}  Let $e \in T_n$ and $f \in T_m$ be evaluable
idempotents
in the generic Temperley Lieb category.  Let $A$ be the $\C(t)$--span in
 $\Hom(m,n)$ of a certain set of $(m,n)$--TL diagrams, and let $A(\tau)$
be
the $\C$--span in  $\Hom(m,n)(\tau)$ of the same set of diagrams.
Then
$$
\dim_{\C(t)} e A f = \dim_\C e(\tau) A(\tau) f(\tau).
$$
\end{proposition}

\begin{proof}  Let $X$ denote the set of TL  diagrams spanning $A$.
Clearly 
$$\dim_{\C(t)} A = \dim_\C A(\tau) = |X|.$$
  Choose a basis
of $e(\tau) A(\tau) f(\tau)$  of the form
$
\{e(\tau)
x f(\tau) : x \in X_0\},
$
where $X_0$ is some subset of $X$.  If the set
$
\{ exf  : x \in X_0\},
$
 were linearly dependent over
$\C(t)$,
then it would be linearly dependent over $\C[t]$, and evaluating at
$\tau$
would give a linear dependence of
$
\{e(\tau)
x f(\tau) : x \in X_0\}
$
over $\C$.
  It
follows
that
$$
\dim_{\C(t)} e A f  \ge \dim_\C e(\tau) A(\tau) f(\tau).
$$
But one has similar inequalities with $e$ replaced by $\one - e$ and/or
$f$
replaced by $\one - f$.  If any of the inequalities were strict, then
adding
them
would give $\dim_{\C(t)} A > \dim_\C A(\tau)$, a contradiction.
\end{proof}

\subsection{The Markov trace} The Markov trace $\Tr=\Tr_n$  is defined
 on $T_n$  (or on $T_n(\tau)$) by the following picture, which
represents an
element in $\End_0 \cong \C(t)$  (resp. $\End(0) \cong \C$).

\begin{figure}[h]

\setlength{\unitlength}{0.00087489in}
\begingroup\makeatletter\ifx\SetFigFont\undefined%
\gdef\SetFigFont#1#2#3#4#5{%
  \reset@font\fontsize{#1}{#2pt}%
  \fontfamily{#3}\fontseries{#4}\fontshape{#5}%
  \selectfont}%
\fi\endgroup%
{\renewcommand{\dashlinestretch}{30}
\begin{picture}(4271,2267)(0,-10)
\path(12,1531)(1137,1531)(1137,631)
    (12,631)(12,1531)
\drawline(1362,1801)(1362,1801)
\thicklines
\path(912,1531)(912,1532)(912,1537)
    (913,1548)(913,1566)(915,1590)
    (916,1617)(918,1645)(921,1671)
    (924,1695)(928,1716)(932,1734)
    (937,1750)(943,1765)(950,1778)
    (956,1790)(964,1801)(973,1812)
    (982,1822)(993,1832)(1004,1842)
    (1017,1851)(1030,1859)(1043,1867)
    (1057,1873)(1070,1879)(1084,1883)
    (1098,1887)(1111,1889)(1124,1891)
    (1137,1891)(1147,1891)(1158,1890)
    (1168,1888)(1179,1885)(1190,1882)
    (1201,1878)(1212,1872)(1223,1866)
    (1234,1858)(1244,1850)(1255,1840)
    (1265,1829)(1274,1818)(1283,1805)
    (1292,1792)(1299,1777)(1307,1762)
    (1313,1746)(1319,1729)(1325,1711)
    (1329,1695)(1332,1678)(1336,1660)
    (1339,1640)(1342,1619)(1344,1595)
    (1347,1570)(1349,1541)(1351,1511)
    (1352,1478)(1354,1442)(1356,1405)
    (1357,1366)(1358,1326)(1359,1287)
    (1360,1251)(1361,1217)(1361,1188)
    (1361,1164)(1362,1147)(1362,1135)
    (1362,1129)(1362,1126)
\path(912,631)(912,630)(912,625)
    (913,614)(914,597)(916,577)
    (918,555)(922,535)(925,517)
    (930,502)(935,488)(942,477)
    (950,466)(957,457)(966,449)
    (977,441)(988,433)(1001,426)
    (1015,419)(1030,412)(1045,407)
    (1061,402)(1076,398)(1092,395)
    (1107,393)(1122,391)(1137,391)
    (1148,391)(1160,392)(1172,394)
    (1184,396)(1196,399)(1208,404)
    (1220,409)(1232,416)(1244,424)
    (1256,433)(1267,443)(1277,455)
    (1287,467)(1296,481)(1304,496)
    (1312,512)(1318,530)(1325,548)
    (1329,564)(1332,580)(1336,597)
    (1339,616)(1342,637)(1344,660)
    (1347,686)(1349,713)(1351,743)
    (1352,776)(1354,811)(1356,849)
    (1357,887)(1358,927)(1359,965)
    (1360,1002)(1361,1035)(1361,1064)
    (1361,1088)(1362,1105)(1362,1117)
    (1362,1123)(1362,1126)
\path(687,631)(687,630)(687,625)
    (688,614)(689,597)(691,576)
    (693,554)(697,533)(700,514)
    (705,498)(710,484)(717,471)
    (725,458)(732,448)(741,438)
    (752,428)(763,417)(776,407)
    (790,397)(805,387)(820,377)
    (836,368)(851,359)(867,351)
    (882,344)(897,337)(912,331)
    (927,325)(942,320)(957,315)
    (974,310)(990,305)(1007,301)
    (1025,298)(1042,294)(1059,292)
    (1075,290)(1092,288)(1107,287)
    (1122,286)(1137,286)(1152,286)
    (1167,287)(1182,288)(1199,290)
    (1215,292)(1232,294)(1250,298)
    (1267,301)(1284,305)(1300,310)
    (1317,315)(1332,320)(1347,325)
    (1362,331)(1373,336)(1385,341)
    (1397,347)(1409,353)(1421,360)
    (1433,368)(1445,377)(1457,387)
    (1469,398)(1481,410)(1492,423)
    (1502,437)(1512,452)(1521,468)
    (1529,485)(1537,502)(1543,521)
    (1550,541)(1554,557)(1557,574)
    (1561,592)(1564,612)(1567,633)
    (1569,657)(1572,682)(1574,711)
    (1576,741)(1577,774)(1579,810)
    (1581,847)(1582,886)(1583,926)
    (1584,965)(1585,1001)(1586,1035)
    (1586,1064)(1586,1088)(1587,1105)
    (1587,1117)(1587,1123)(1587,1126)
\path(462,631)(462,630)(462,627)
    (463,619)(465,607)(467,590)
    (471,570)(475,550)(480,530)
    (486,512)(493,495)(502,480)
    (512,465)(523,451)(537,436)
    (548,426)(560,415)(573,404)
    (587,393)(603,381)(620,369)
    (638,357)(657,345)(677,332)
    (697,320)(719,308)(741,297)
    (762,286)(784,275)(806,265)
    (828,256)(849,247)(871,240)
    (891,232)(912,226)(933,220)
    (954,215)(975,210)(997,206)
    (1019,202)(1042,199)(1065,197)
    (1089,195)(1113,194)(1137,193)
    (1161,194)(1185,195)(1209,197)
    (1232,199)(1255,202)(1277,206)
    (1299,210)(1320,215)(1341,220)
    (1362,226)(1383,233)(1403,240)
    (1425,248)(1446,257)(1468,267)
    (1490,278)(1512,290)(1533,303)
    (1555,317)(1577,332)(1597,347)
    (1617,364)(1636,382)(1654,400)
    (1671,418)(1687,437)(1701,457)
    (1714,477)(1726,497)(1737,518)
    (1745,537)(1753,556)(1760,576)
    (1766,598)(1771,621)(1777,646)
    (1781,673)(1786,703)(1789,734)
    (1793,769)(1796,805)(1799,843)
    (1802,883)(1804,923)(1806,963)
    (1808,1000)(1809,1034)(1810,1064)
    (1811,1087)(1811,1105)(1812,1117)
    (1812,1123)(1812,1126)
\path(237,631)(237,630)(238,627)
    (239,619)(241,605)(245,587)
    (250,567)(256,545)(264,523)
    (273,503)(284,483)(297,464)
    (311,445)(329,426)(350,406)
    (363,394)(377,381)(393,368)
    (410,354)(428,340)(448,325)
    (468,310)(490,294)(512,278)
    (536,262)(560,246)(585,230)
    (610,214)(636,199)(661,183)
    (687,169)(712,155)(737,141)
    (761,128)(785,116)(808,105)
    (831,95)(853,85)(875,76)
    (898,67)(921,59)(944,51)
    (967,45)(991,39)(1014,34)
    (1039,30)(1063,26)(1087,24)
    (1112,22)(1137,22)(1162,22)
    (1187,24)(1211,26)(1235,30)
    (1260,34)(1283,39)(1307,45)
    (1330,51)(1353,59)(1376,67)
    (1400,76)(1421,85)(1443,95)
    (1466,105)(1489,116)(1513,129)
    (1537,142)(1562,156)(1587,170)
    (1613,186)(1638,202)(1664,218)
    (1689,235)(1714,253)(1738,271)
    (1762,288)(1784,306)(1806,324)
    (1826,342)(1846,360)(1864,377)
    (1881,394)(1897,411)(1911,427)
    (1925,443)(1939,463)(1952,483)
    (1964,503)(1975,524)(1984,545)
    (1992,569)(2000,593)(2006,620)
    (2012,648)(2018,678)(2022,709)
    (2026,739)(2030,769)(2032,795)
    (2034,817)(2036,835)(2036,846)
    (2037,853)(2037,856)
\path(687,1531)(687,1532)(687,1535)
    (688,1543)(688,1556)(690,1575)
    (692,1598)(694,1622)(697,1647)
    (701,1671)(705,1694)(710,1714)
    (715,1732)(722,1749)(729,1764)
    (737,1779)(747,1793)(757,1806)
    (767,1819)(779,1831)(792,1844)
    (806,1857)(822,1870)(838,1882)
    (855,1895)(872,1907)(890,1919)
    (908,1930)(926,1941)(944,1951)
    (961,1960)(978,1968)(994,1976)
    (1009,1982)(1025,1988)(1041,1995)
    (1057,2000)(1074,2005)(1091,2008)
    (1107,2012)(1124,2014)(1141,2015)
    (1158,2016)(1175,2015)(1192,2014)
    (1209,2012)(1225,2008)(1241,2005)
    (1256,2000)(1272,1995)(1287,1988)
    (1301,1982)(1315,1976)(1330,1968)
    (1345,1959)(1361,1949)(1376,1938)
    (1392,1926)(1408,1913)(1424,1898)
    (1440,1882)(1455,1865)(1469,1848)
    (1482,1829)(1495,1810)(1506,1790)
    (1517,1770)(1526,1748)(1535,1726)
    (1540,1708)(1546,1690)(1550,1670)
    (1555,1649)(1559,1627)(1562,1602)
    (1565,1576)(1568,1547)(1571,1515)
    (1574,1481)(1576,1445)(1578,1407)
    (1580,1368)(1581,1328)(1583,1289)
    (1584,1251)(1585,1218)(1586,1188)
    (1586,1165)(1587,1147)(1587,1135)
    (1587,1129)(1587,1126)
\path(462,1531)(462,1532)(462,1535)
    (463,1543)(464,1557)(466,1576)
    (468,1600)(471,1625)(475,1651)
    (479,1677)(485,1700)(491,1722)
    (497,1743)(505,1762)(514,1780)
    (525,1798)(537,1816)(548,1831)
    (560,1845)(573,1861)(587,1876)
    (603,1892)(620,1909)(638,1925)
    (657,1942)(677,1958)(697,1975)
    (719,1991)(741,2006)(762,2021)
    (784,2035)(806,2049)(828,2061)
    (849,2072)(871,2083)(891,2092)
    (912,2101)(933,2109)(954,2116)
    (975,2122)(997,2128)(1019,2133)
    (1042,2137)(1065,2140)(1089,2142)
    (1113,2144)(1137,2144)(1161,2144)
    (1185,2142)(1209,2140)(1232,2137)
    (1255,2133)(1277,2128)(1299,2122)
    (1320,2116)(1341,2109)(1362,2101)
    (1381,2093)(1400,2085)(1419,2075)
    (1438,2065)(1458,2053)(1478,2041)
    (1498,2028)(1518,2014)(1537,1998)
    (1557,1982)(1577,1965)(1596,1948)
    (1614,1929)(1631,1911)(1648,1891)
    (1664,1872)(1678,1852)(1692,1832)
    (1705,1811)(1716,1791)(1727,1770)
    (1737,1748)(1745,1728)(1753,1708)
    (1760,1686)(1766,1663)(1771,1639)
    (1777,1613)(1781,1585)(1786,1555)
    (1789,1522)(1793,1487)(1796,1450)
    (1799,1411)(1802,1371)(1804,1330)
    (1806,1290)(1808,1253)(1809,1218)
    (1810,1189)(1811,1165)(1811,1147)
    (1812,1135)(1812,1129)(1812,1126)
\path(237,1531)(237,1532)(237,1535)
    (238,1543)(240,1557)(242,1577)
    (246,1601)(251,1627)(256,1654)
    (263,1680)(271,1705)(280,1728)
    (290,1750)(302,1771)(316,1791)
    (331,1811)(350,1831)(363,1845)
    (377,1859)(393,1874)(410,1889)
    (428,1904)(448,1920)(468,1936)
    (490,1953)(512,1970)(536,1986)
    (560,2003)(585,2020)(610,2036)
    (636,2052)(661,2068)(687,2083)
    (712,2097)(737,2111)(761,2124)
    (785,2136)(808,2147)(831,2157)
    (853,2167)(875,2176)(898,2185)
    (921,2193)(944,2201)(967,2207)
    (991,2213)(1014,2218)(1039,2222)
    (1063,2226)(1087,2228)(1112,2230)
    (1137,2230)(1162,2230)(1187,2228)
    (1211,2226)(1235,2222)(1260,2218)
    (1283,2213)(1307,2207)(1330,2201)
    (1353,2193)(1376,2185)(1400,2176)
    (1421,2167)(1443,2157)(1466,2147)
    (1489,2136)(1513,2124)(1537,2111)
    (1562,2097)(1587,2083)(1613,2068)
    (1638,2053)(1664,2037)(1689,2021)
    (1714,2004)(1738,1988)(1762,1972)
    (1784,1955)(1806,1939)(1826,1924)
    (1846,1908)(1864,1893)(1881,1879)
    (1897,1865)(1911,1852)(1925,1838)
    (1939,1823)(1953,1808)(1965,1793)
    (1976,1778)(1986,1763)(1995,1748)
    (2003,1732)(2010,1716)(2016,1699)
    (2021,1682)(2026,1665)(2029,1647)
    (2032,1629)(2034,1611)(2035,1593)
    (2036,1574)(2037,1555)(2037,1535)
    (2037,1515)(2037,1493)(2037,1475)
    (2037,1456)(2037,1435)(2037,1413)
    (2037,1389)(2037,1363)(2037,1335)
    (2037,1304)(2037,1271)(2037,1235)
    (2037,1196)(2037,1156)(2037,1114)
    (2037,1071)(2037,1029)(2037,990)
    (2037,954)(2037,922)(2037,897)
    (2037,878)(2037,866)(2037,859)(2037,856)
\put(462,1126){\makebox(0,0)[lb]
{\smash{{{\SetFigFont{12}{14.4}{\rmdefault}{\mddefault}{\updefault}$a$}}
}}}
\put(2487,1171){\makebox(0,0)[lb]
{\smash{{{\SetFigFont{12}{14.4}{\rmdefault}{\mddefault}{\updefault}$ =
\Tr(a) \in 
 \End(0) $}}}}}
\end{picture}
}

\caption{The categorical trace of an element $a \in T_n$.}

\end{figure}
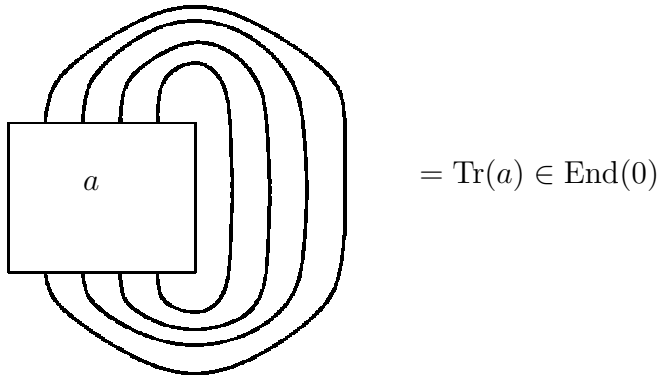

 On an $(n,n)$--TL diagram
$a\in T_n$, the trace is evaluated
 geometrically by closing  up the diagram as in the figure, and counting
the
number
$c(a)$ of components (closed loops); then $\Tr(a)=d^{c(a)}$.

It will be useful to give the following inductive description
of closing up a diagram. We define a map $\ve_n:T_{n+1}\to T_n$
(known as a conditional expectation in operator algebras) by
only closing up the last strand; algebraically it can be defined by
$$a\in T_{n+1}\quad \mapsto \quad (1_n\otimes \cup)\circ(a\otimes 1)
\circ (1_n\otimes \cap).$$
\ignore{where $1_i$ indicates $i$ vertical strands, and composition of
our morphism from right to left corresponds to putting the graphs
on top of each other.} If $k>n$, the map $\ve_{n,k}$ is defined by
$\ve_{n,k}=\ve_{n}\circ\ve_{n+1}\ ...\ \circ\ve_{k-1}$.
It follows from the definitions that $\Tr(a)=\ve_{0,n}$ for $a\in T_n$.

It is well-known that $\Tr$ is indeed a functional
satisfying $\Tr(ab)=\Tr(ba)$; one easily checks that this equality is
even
true if $a\in \Hom(n,m)$ and $b\in \Hom(m,n)$. We need the following
well-known fact:

\begin{lemma}\label{condexp}
Let $f\in T_{n+m}$ and let $p\in T_n$ such that $(p\otimes
1_m)f(p\otimes
1_m)
=f$, where $p$ is a minimal idempotent in $T_n$.
Then $\ve_{n,n+m}(f)=\gamma p$, where $\gamma = Tr_{n+m}(f)/Tr_n(p)$
\end{lemma}

\begin{proof}
It follows from the definitions that
$$p\ve_{n,n+m}(f)p=\ve_{n,n+m}((p\otimes 1_m)f(p\otimes
1_m))=\ve_{n,n+m}(f).$$
As $p$ is a minimal idempotent in $T_n$, $\ve_{n,n+m}(f)=\gamma p,$  for
some scalar $\gamma$. Moreover, by our definition of trace, we have
$Tr_{n+m}(f) =Tr_n(\ve_{n,n+m}(f))=\gamma Tr_n(p)$. This determines the
value of $\gamma$.
\end{proof}

The negligible morphisms
$\Ne(n,m)$ are defined to be all elements $a\in \Hom(n,m)$
for which $Tr(ab)=0$ for all $b\in \Hom(m,n)$. It is well-known
that the negligible morphisms form an ideal in $\TL$.

\section{The structure of the Temperley Lieb algebras}

\subsection{The generic Temperley Lieb algebras}

Recall that a {\em Young  diagram} $\lambda = [\lambda_1,\lambda_2,\
...\
\lambda_k]$ is a left justified array of boxes with $\lambda_i$ boxes in
the $i$-th row and $\lambda_i\geq \lambda_{i+1}$ for all $i$.
For example,
$$
[5,3] = {\YD \R5,\R3,.}.
$$
{\em All Young diagrams in this paper will have at most two rows.}
For $\la$ a Young diagram with $n$ boxes, a
{\em  Young tableau  of shape $\la$} is
a filling of
$\lambda$ with the numbers 1 through $n$  so that
the numbers increase in each row and column.   The number of Young
tableax
of
shape $\la$ is denoted by $f_\la$.

The generic Temperley Lieb algebras $T_n$ are known (~\cite{jones}) to
decompose as
direct sums of full matrix algebras over the field $\C(t)$,
$T_n = \bigoplus_\la T_\la$, where the sum is over all  Young diagrams
$\lambda$
with $n$ boxes (and with no more than two rows), and $T_\lambda$ is
isomorphic
to an
$f_\la$-by-$f_\la$ matrix algebra.

When $\lambda$ and $\mu$  are Young diagrams of size $n$ and $n+1$,
one has a
(non-unital) homomorphism of $T_\lambda$ into $T_\mu$ given by
$x \mapsto (x \otimes 1)z_\mu$, where $z_\mu$ denotes the minimal
central
idempotent in $T_{n+1}$ such that $T_\mu = T_{n+1} z_\mu$.  Let
$\glm$ denote the rank of
$ (e \otimes 1)z_\mu$, where $e$ is any minimal idempotent in
$T_\lambda$.  It is known that $\glm = 1$ in case $\mu$ is obtained from
$\la$ by adding one box, and $\glm = 0$ otherwise.

\begin{figure}
\begin{center}
\setlength{\unitlength}{0.0125in}
\begin{picture}(348,304)(0,-10)
\drawline(0,280)(280,0)
\drawline(120,160)(0,40)(40,0)(160,120)
\drawline(80,200)(0,120)(120,0)(200,80)
\drawline(40,240)(0,200)(200,0)(240,40)
\put(240,40){\makebox(0,0)[lb]{\raisebox{0pt}[0pt][0pt]
{\shortstack[l]{{\twlrm $\YD \R6,.$}}}}}
\put(50,0){\makebox(0,0)[lb]{\raisebox{0pt}[0pt][0pt]
{\shortstack[l]{{\twlrm $\YD \R4,\R3,.$}}}}}
\put(280,0){\makebox(0,0)[lb]{\raisebox{0pt}[0pt][0pt]
{\shortstack[l]{{\twlrm $\YD \R7,.$}}}}}
\put(215,0){\makebox(0,0)[lb]{\raisebox{0pt}[0pt][0pt]
{\shortstack[l]{{\twlrm $\YD \R6,\R1,.$}}}}}
\put(130,0){\makebox(0,0)[lb]{\raisebox{0pt}[0pt][0pt]
{\shortstack[l]{{\twlrm $\YD \R5,\R2,.$}}}}}
\put(10,40){\makebox(0,0)[lb]{\raisebox{0pt}[0pt][0pt]
{\shortstack[l]{{\twlrm $\YD \R3,\R3,.$}}}}}
\put(175,40){\makebox(0,0)[lb]{\raisebox{0pt}[0pt][0pt]
{\shortstack[l]{{\twlrm $\YD \R5,\R1,.$}}}}}
\put(95,40){\makebox(0,0)[lb]{\raisebox{0pt}[0pt][0pt]
{\shortstack[l]{{\twlrm $\YD \R4,\R2,.$}}}}}
\put(130,75){\makebox(0,0)[lb]{\raisebox{0pt}[0pt][0pt]
{\shortstack[l]{{\twlrm $\YD \R4,\R1,.$}}}}}
\put(50,75){\makebox(0,0)[lb]{\raisebox{0pt}[0pt][0pt]
{\shortstack[l]{{\twlrm $\YD \R3,\R2,.$}}}}}
\put(10,115){\makebox(0,0)[lb]{\raisebox{0pt}[0pt][0pt]
{\shortstack[l]{{\twlrm $\YD \R2,\R2,.$}}}}}
\put(90,115){\makebox(0,0)[lb]{\raisebox{0pt}[0pt][0pt]
{\shortstack[l]{{\twlrm $\YD \R3,\R1,.$}}}}}
\put(50,155){\makebox(0,0)[lb]{\raisebox{0pt}[0pt][0pt]
{\shortstack[l]{{\twlrm $\YD \R2,\R1,.$}}}}}
\put(205,80){\makebox(0,0)[lb]{\raisebox{0pt}[0pt][0pt]
{\shortstack[l]{{\twlrm $\YD \R5,.$}}}}}
\put(165,115){\makebox(0,0)[lb]{\raisebox{0pt}[0pt][0pt]
{\shortstack[l]{{\twlrm $\YD \R4,.$}}}}}
\put(125,155){\makebox(0,0)[lb]{\raisebox{0pt}[0pt][0pt]
{\shortstack[l]{{\twlrm $\YD \R3,.$}}}}}
\put(10,275){\makebox(0,0)[lb]{\raisebox{0pt}[0pt][0pt]
{\shortstack[l]{{\twlrm $\emptyset$}}}}}
\put(45,240){\makebox(0,0)[lb]{\raisebox{0pt}[0pt][0pt]
{\shortstack[l]{{\twlrm $\YD \R1,.$}}}}}
\put(85,200){\makebox(0,0)[lb]{\raisebox{0pt}[0pt][0pt]
{\shortstack[l]{{\twlrm $\YD \R2,.$}}}}}
\put(10,200){\makebox(0,0)[lb]{\raisebox{0pt}[0pt][0pt]
{\shortstack[l]{{\twlrm $\YD \R1,\R1,.$}}}}}
\end{picture}

\caption{Bratteli diagram for the sequence $(T_n)$}
\label{bratteli}
\end{center}
\end{figure}
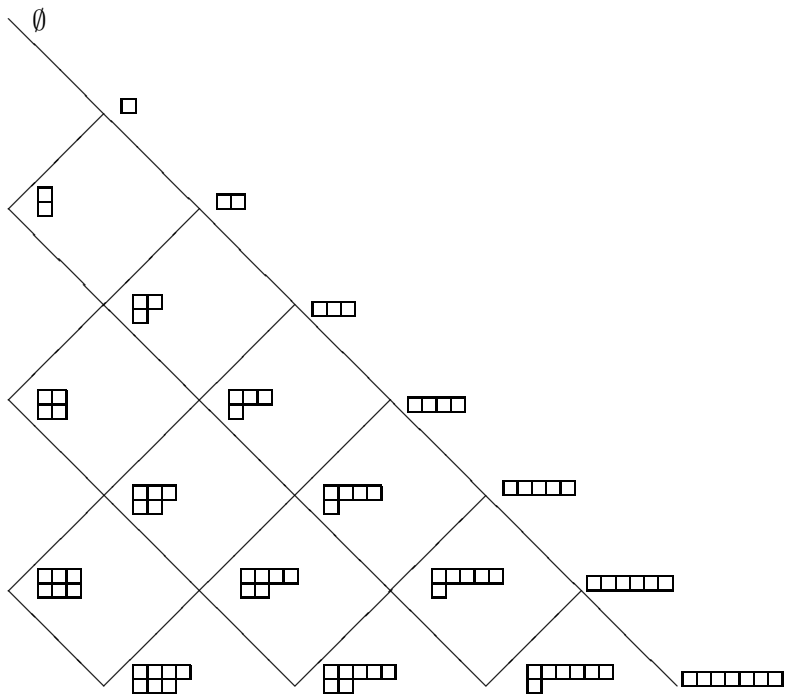

One can describe the embedding of $T_n$ into $T_{n+1}$ by a {\em
Bratteli
diagram}  (or induction-restriction diagram), which is a bipartite graph
with
vertices labelled by two-row Young diagrams of size $n$ and $n+1$
(corresponding to the simple components of
$T_n$ and $T_{n+1}$) and with $\glm$ edges joining the vertices labelled
by
$\la$ and $\mu$.  That is $\lambda$ and $\mu$ are  joined by an edge
precisely when $\mu$ is obtained from $\la$ by adding one box.  The
sequence
of embeddings $T_0 \rightarrow T_1 \rightarrow T_2 \rightarrow \dots$ is
described by a multilevel Bratteli diagram, as shown in Figure
\ref{bratteli}.

 A  tableau of shape $\la$
may be identified with an
increasing sequence of Young diagrams beginning with the empty diagram
and
ending at $\la$; namely  the $j$-th  diagram in the sequence is the
subdiagram of $\la$ containing the numbers 1, 2, \dots, $j$. Such a
sequence
may
also be interpreted as a {\em path} on the Bratteli diagram of Figure
\ref{bratteli}, beginning at  the empty diagram and ending at
$\la$.

\subsection{Path idempotents}
One can define a familiy of minimal idempotents $p_t$ in $T_n$, labelled
by
paths $t$  of length $n$ on the Bratteli diagram (or equivalently, by
Young
tableaux of size $n$), with the following
properties:

\begin{enumerate}
\item $p_t p_s = 0$ if $t, s$ are different paths both of length $n$.
\item $z_\la = \sum \{p_t : t \text{ ends at } \la\}$.
\item $p_t \otimes 1 = \sum \{p_s : s  \text{ has length $n+1$ and
extends } t \}$
\end{enumerate} 

Let $t$ be a path of length $n$ and shape $\lambda$ and let
$\mu$  be a Young diagram of size $n + m$.  It follows that $(p_t
\otimes
1_m) z_\mu \ne 0$ precisely when there is a path on the Bratteli diagram
from $\la$ to $\mu$.
It has been shown in ~\cite{jones} that (in our notations)
$\Tr(p_t)=[\la_1-\la_2+1]$, where $[m]=(t^m-t^{-m})/(t-t^{-1})$
for any integer $m$, and where $\la$ is the endpoint of the path $t$.
Observe that we get the same value for diagrams $\la$ and $\mu$
which are in the same column in the Bratteli diagram.

The idempotents $p_t$ were defined by recursive formulas in ~\cite{type
A},
generalizing the formulas for the Jones-Wenzl idempotents in
~\cite{sequence
of projections}.  

\subsection{Specializations at  non-roots of unity}

When $\tau$ is not a proper root of unity,  the Temperley Lieb algebras
$T_n(\tau)$ are semi-simple complex algebras with the ``same" structure
as
generic Temperley Lieb algebras.  That is,
$T_n(\tau) = \bigoplus_\la T_\la(\tau)$, where  $T_\lambda(\tau)$ is
isomorphic to an
$f_\la$-by-$f_\la$ matrix algebra over $\C$.  The embeddings
$T_n(\tau) \rightarrow T_{n+1}(\tau)$ are described by the Bratteli
diagram
as before.  The idempotents $p_t$, and the minimal central idempotents
$z_\lambda$, in the generic algebras $T_n$,  are
evaluable at
$\tau$, and the  evaluations $p_t(\tau)$, resp. $z_\lambda(\tau)$,
satisfy
analogous properties.

\subsection{Specializations at roots of unity and evaluable idempotents}
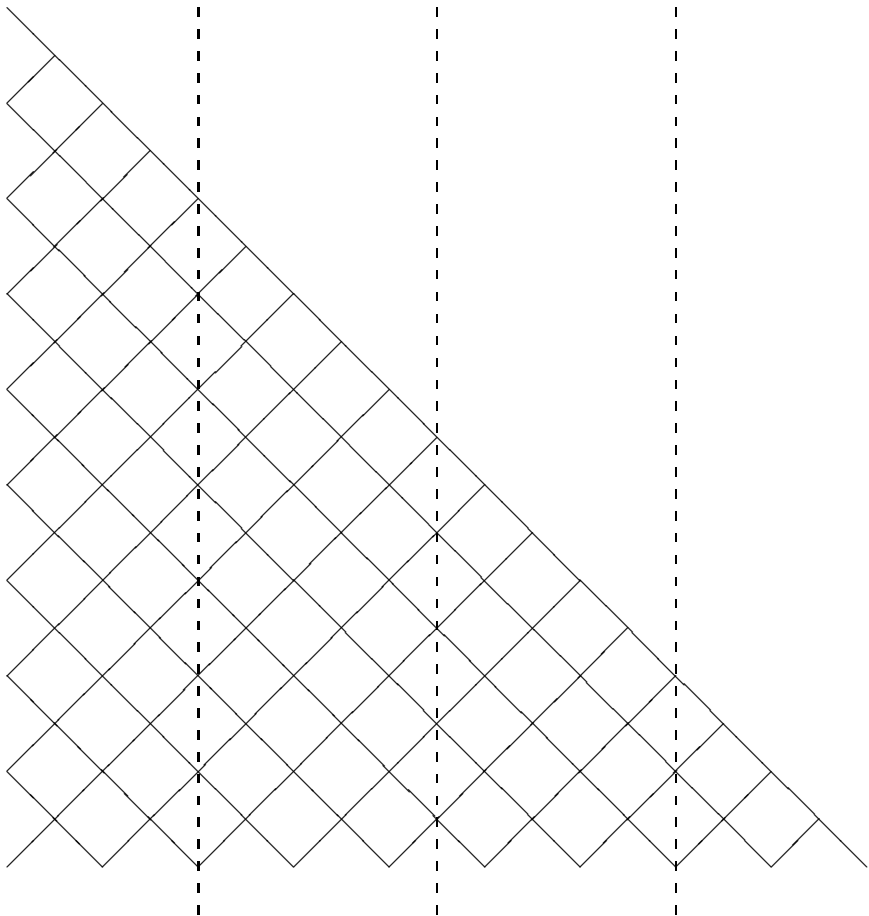
\begin{figure}[ht]
\begin{center}
\setlength{\unitlength}{0.0125in}
\begin{picture}(360,395)(0,-10)
\dashline{4.000}(280,380)(280,0)
\dashline{4.000}(180,380)(180,0)
\dashline{4.000}(80,380)(80,0)
\drawline(340,40)(320,20)
\drawline(320,60)(280,20)
\drawline(300,80)(240,20)
\drawline(280,100)(200,20)
\drawline(260,120)(160,20)
\drawline(240,140)(120,20)
\drawline(220,160)(80,20)
\drawline(200,180)(40,20)
\drawline(180,200)(0,20)
\drawline(160,220)(0,60)(40,20)
\drawline(140,240)(0,100)(80,20)
\drawline(120,260)(0,140)(120,20)
\drawline(100,280)(0,180)(160,20)
\drawline(80,300)(0,220)(200,20)
\drawline(60,320)(0,260)(240,20)
\drawline(60,320)(60,320)
\drawline(40,340)(0,300)(280,20)
\drawline(20,360)(0,340)(320,20)
\drawline(0,380)(360,20)
\end{picture}
\end{center}
\caption{Critical lines}
\label{critical lines}
\end{figure}

We require some terminology for discussing the case where $\tau$ is a
root of unity. Let $q = \tau^2$, and suppose that $q$ is a primitive
$\ell$-th root of unity.  We say that a Young diagram $\la$  is
{\em critical} if $w(\la) :=\la_1 - \la_2 +1$ is divisible by $\ell$.
The $m$-th {\em critical line} on the Bratteli diagram for the generic
Temperly
Lieb algebra is the line containing the diagrams $\la$ with $w(\la) =
ml$.
See
Figure \ref{critical lines}.

Say that two non-critical diagrams $\la$ and
$\mu$
with the same number of boxes are {\em reflections of one another in the
$m$-th
critical line} if
$\la \ne \mu$ and $|w(\la) - m\ell| = |w(\mu) - m\ell| < \ell$.
(For example, with $\ell = 3$, $[2,2]$ and $[4]$ are reflections in the
first
critical line $w(\la) = 3$.)

For $\tau$ a proper root of unity, the formulas for path idempotents in
\cite{sequence  of projections} and \cite{type A} generally contain
poles
at $\tau$, i.e. the idempotents are not evaluable. However, suitable
sums
of path idempotents are evaluable.

Suppose $w(\la) \le \ell$ and $t$ is a path of shape $\la$
which stays strictly
to the left of the first critical line (in case $w(\la) < \ell$), or
hits
the
first critical line for the first time at $\la$ (in case $w(\la) =
\ell$);
then
$p_t$ is evaluable at
$\tau$, and  furthermore $\Tr(p_t) = [w(\la)]_\tau =
(\tau^{w(\la)}-\tau^{-{w(\la)}})/(\tau-\tau^{-1})$.

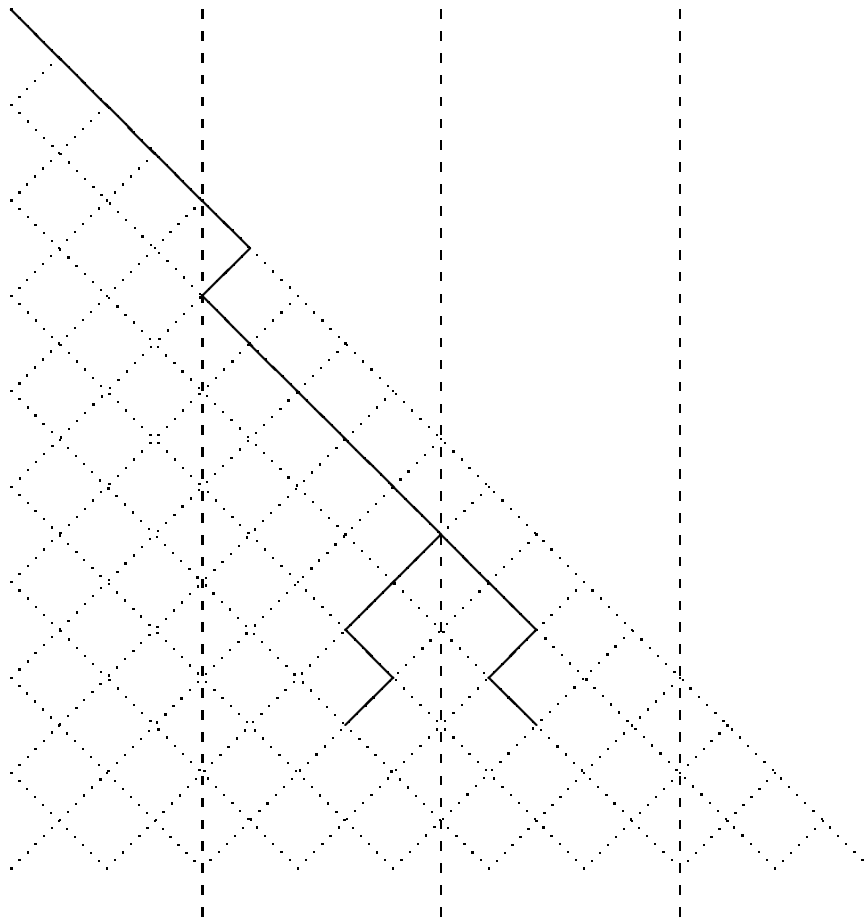
\begin{figure}[ht]
\setlength{\unitlength}{0.0125in}
\begin{picture}(360,395)(0,-10)
\dottedline{5}(40,340)(0,300)(280,20)
\thicklines
\drawline(0,380)(80,300)(100,280)
    (80,260)(180,160)(140,120)
    (160,100)(140,80)
\drawline(180,160)(220,120)(200,100)(220,80)
\thinlines
\dottedline{5}(320,20)(0,340)(20,360)
\dottedline{5}(360,20)(0,380)
\dottedline{5}(320,60)(280,20)
\dottedline{5}(340,40)(320,20)
\dottedline{5}(280,100)(200,20)
\dottedline{5}(260,120)(160,20)
\dottedline{5}(300,80)(240,20)
\dottedline{5}(220,160)(80,20)
\dottedline{5}(200,180)(40,20)
\dottedline{5}(240,140)(120,20)
\dottedline{5}(160,220)(0,60)(40,20)
\dottedline{5}(140,240)(0,100)(80,20)
\dottedline{5}(120,260)(0,140)(120,20)
\dottedline{5}(100,280)(0,180)(160,20)
\dottedline{5}(80,300)(0,220)(200,20)
\dottedline{5}(60,320)(0,260)(240,20)
\dottedline{5}(180,200)(0,20)
\dashline{4.000}(80,380)(80,0)
\dashline{4.000}(180,380)(180,0)
\dashline{4.000}(280,380)(280,0)
\drawline(60,320)(60,320)
\end{picture}
\caption{A path and its reflected path.}
\label{conjugate paths}
\end{figure}

 For each
critical diagram $\la$  of size $n$, the minimal central idempotent
$z_\la$
in 
$T_n$ is evaluable at $\tau$.  Furthermore, for each non-critical
diagram
$\la$ of size $n$, an evaluable idempotent $z_\la^L =\sum p_{t}\in T_n$
was defined in ~\cite{pacific} as follows: The summation goes over all
paths $t$ ending in $\la$ for which the last critical line hit by $t$ is
the
one nearest to $\la$ to the left {\it and} over the paths obtained from
such
$t$ by reflecting its part after the last critical line in the critical
line
(see Figure \ref{conjugate paths}).

These idempotents have the
following properties (which were shown in \break ~\cite{pacific}):
\begin{enumerate}
\item $\{z_\la(\tau) : \la \text{ critical }\} \cup
\{z_\mu^L(\tau) : \mu \text{ non-critical }\}$ is a partition
of unity in 
$T_n(\tau)$;  that is, the idempotents are mutually orthogonal and sum
to
the identity.
\item   $z_\la(\tau)$  is a  minimal central idempotent in $T_n(\tau)$
if $\la$ is critical, and
$z_\la^L(\tau)$ is minimal
central modulo the nilradical of $T_n$ if $\la$ is not critical
  (see
~\cite{pacific}, Theorem 2.2 and Theorem 2.3).
\item  For $\la$ and $\mu$ non-critical, $z_\la^L(\tau) T_n(\tau)
z_\mu^L(\tau) \ne 0$ only if $\la =\mu$, or if there is exactly
one critical line between
$\la$ and  $\mu$ which reflects $\la$ to $\mu$.
If in this case $\mu$ is to the left of $\la$,
$z_\la^L T_n z_{\la'}^L \subseteq T_\mu$
 (in the generic Temperley Lieb algebra).
\item Let $z_n^{reg}=\sum p_t$, where the summation goes over all paths
$t$ which stay strictly to the left of the first critical line,
and let $z_n^{nil}=\one-z_n^{reg}$.
Then both $z_n^{reg}$ and $z_n^{nil}$ are evaluable; this is a direct
consequence of the fact that $z_n^{reg}=\sum_\la z_\la^L$,
where the summation goes over diagrams $\la$ with $n$ boxes with
width $w(\la)<\ell$.
\end{enumerate}
\ignore{
Finally, we also mention that the restriction rule in the root of unity
case for the semisimple quotients of $T_n$ is very similar to the
generic case. It coincides in all cases except if one of the
involved diagrams lies on a critical line. If a predecessor of
$\la$ lies in a critical line
a simple $T_{n,\la}$-module $W$ decomposes into the direct sum of
3 simple $T_{n-1}$ modules.}

\begin{proposition}
The ideal of negligible morphisms in $\TL(\tau)$ is generated by the
idempotent $p_{[\ell-1]}(\tau) \in T_{\ell-1}(\tau)$.
\end{proposition}

\begin{proof} Let us first show
that $z_n^{nil}(\tau)$ is in the ideal generated by
 $p_{[\ell-1]}(\tau)$ for all $n$. This is clear for $n<\ell$,
as $z_{\ell -1}^{nil}=p_{[\ell-1]}$ and $z_n^{nil}=0$ for $n<\ell-1$.

Moreover, $z^{nil}_n$ is a central idempotent in the maximum semisimple
quotient of $T_n$, whose minimal central idempotents are the
$z_\la^L$ with $w(\la)\geq \ell$. One checks pictorially that
$p_{[\ell-1]}z_\la^L\neq 0$ for any such $\lambda$ (i.e.
the path to $[\ell-1]$ can be extended to a path $t$ for which $p_t$
is a summand of $z_\la^L$).  This proves our assertion
in the maximum semisimple quotient of $T_n$; it is well-known that
in this case also the idempotent itself must be in the ideal generated
by $p_{[\ell-1]}$. In particular,  $\Hom(n,m)z_m^{nil}(\tau)
+z_n^{nil}(\tau)\Hom(n,m)$ is also contained in this ideal.

By [GW], Theorem 2.2 (c), for $\la$ a Young diagram of size $n$,
with $w(\la) < \ell$,
$z_\la^L T_n z_\la^L(\tau)$ is a full matrix algebra, which moreover
contains a minimal idempotent $p_t$ of trace $\Tr(p_t) =
[w(\la)]_\tau \ne 0$.  Therefore $$z_\la^L T_n z_\la^L(\tau) \cap
\Neg(n,n) = (0).$$   Furthermore, $z_n^{reg} T_n z_n^{reg}(\tau) =
\sum z_\la^L T_n  z_\la^L(\tau) $, by Fact 4 above, so
$$z_n^{reg} T_n z_n^{reg}(\tau) \cap \Neg(n,n) = (0)$$
as well.  Now for $x \in \Neg(n,n)$,
one has $z_n^{reg}(\tau) x z_n^{reg}(\tau) = 0$, so
$$x \in T_n(\tau) z_n^{nil}(\tau)+z_n^{nil}(\tau)T_n (\tau).$$

We have shown that $\Neg(n,n)$ is contained in the ideal of $\TL(\tau)$
generated
by $p_{[l-1]}$, for all $n$.    That the same is true for $\Neg(m,n)$
with $n\neq m$ follows
from using the embeddings, and their left inverses, described at the end
of
Section 1.1.

\ignore{
On the other hand, if $n=m$, its complement $z_n^{reg}T_nz_n^{reg}$
(as a linear direct summand) has the same dimension as the quotient
of $T_n/\Ne(n,n)$: this can be read off inductively from the
Bratteli diagram of the quotient and the definition of $z_n^{reg}$.
Hence $T_nz_n^{nil}+z_n^{nil}T_n$ must contain all negligible morphisms
in $T_n$. The claim for the general case follows from this and
the embedding described at the end of Section 1.1
}

\end{proof}

\section{Ideals}

\begin{proposition} Any proper ideal in $\TL$ (or in $\TL(\tau)$) is
contained in the ideal of negligible morphisms.
\end{proposition}

\begin{proof} Let $a \in \Hom(m, n)$.  For all $b \in \Hom(n,m)$,
$\tr(ba)$ is in the intersection of the ideal generated by $a$ with the
scalars $\End(0)$.  If $a$ is not negligible, then
the ideal generated by $a$ contains an non-zero scalar, and therefore
contains all morphisms.
\end{proof}

\begin{corollary}  The categories $\TL$ and $\TL(\tau)$ for $\tau$ not a
proper root of unity have no non-zero proper ideals.
\end{corollary}

\begin{proof}  There are no non-zero negligible morphisms in
$\TL$ and in $\TL(\tau)$ for $\tau$ not a
proper root of unity.
\end{proof}

\ignore{
\begin{lemma} \label{embedding lemma}
For  $m, n \in \nat_0$,
there is an injective linear map $$\iota : \Hom(m,n)(\tau) \rightarrow
T_{\max(m,n)}(\tau)$$ such that for any ideal
$J$ in $\TL(\tau)$,  $$\iota(J \cap \Hom(m,n)) \subseteq J \cap
T_{\max(m,n)}(\tau).$$
\end{lemma}
\begin{proof}  blah, blah
\end{proof}
}

\begin{theorem}  Suppose that $\tau$ is a proper root of unity.  Then
the negligible morphisms form the unique non-zero proper ideal in
$\TL(\tau)$.
\end{theorem}

\begin{proof}  Let $J$ be a non-zero proper ideal in $\TL(\tau)$.
By the embeddings discussed at the end of Section 1.1, we can assume
$J\cap T_n\neq 0$ for some $n$.

Now let $a$ be a non-zero element of $J \cap T_n(\tau)$.  Since
$\{z_\la(\tau)  \}
\cup
\{z_\mu^L(\tau) \}$  is a partition
of unity in $T_n(\tau)$, one of the following conditions hold:
{\renewcommand{\labelenumi}{{(\alph{enumi})}}
\begin{enumerate}
\item $b = a z_\mu(\tau) \ne 0$ for some
critical diagram
$\mu$.
\item  
$b = z_\mu^L(\tau) a z_\mu^L(\tau) \ne 0$ for some non-critical diagram
$\mu$.
\item  
$b = z_\la^L(\tau) a z_{\la'}^L(\tau) \ne 0$
 for some pair $\la, \la'$ of non-critical diagrams which are
reflections of
one another in a critical line.  In this  case, let $\mu$ denote the
leftmost
of the two diagrams $\la, \la'$.
\end{enumerate}
}

In each of the three cases, one has $b \in e(\tau) T_n(\tau) f(\tau)$,
where
$e, f$ are evaluable idempotents in $T_n$ such that $e T_n f  \subseteq
T_\mu$.
Let $\alpha$ be a Young diagram on the
first critical line  of size $n+m$, such that there exists a path on the
generic
Bratteli diagram connecting $\mu$ and $\alpha$.  Then one has
\begin{align*}
\dim_\C\ &z_\alpha(\tau)(e(\tau) \otimes 1_m) ( T_n(\tau)  \otimes \C\,
1_m)
(f(\tau)
\otimes 1_m) \\ &=
\dim_{\C(t)}z_\alpha (e \otimes \id_m )( T_n  \otimes  \C(t) 1_m)
(f \otimes 1_m )
\\
&= \dim_{\C(t)} e T_n f =
\dim_\C e(\tau) T_n(\tau) f(\tau)
\end{align*}
where the first and last equalities result from the principle of
constancy
of
dimension, and the second equality is because $x \mapsto  z_\alpha(x
\otimes
1_m) $ is injective from $T_\mu$ to $T_\alpha$.
But then it follows that $ x \mapsto z_\alpha(\tau) (x \otimes 1_m)$ is
injective
on $e(\tau) T_n(\tau) f(\tau)$.  In particular
$(b \otimes 1_m)z_\alpha $ is a non-zero element of $J \cap T_\alpha$.
Hence there exists $c\in T_\al$ such that $f=c(b\otimes 1_m)z_\al$
is an idempotent. After conjugating (and multiplying with $p_{[\ell-1]}
\otimes 1_m$, if necessary), we can assume $f$ to be a subidempotent
of  $p_{[\ell-1]} \otimes 1_m$. But then $\ve_{\ell-1+m,\ell-1}(f)$
is a multiple of  $p_{[\ell-1]}$, by Lemma \ref{condexp},
with the multiple equal to the rank of $f$ in $T_\al$.
This, together with Prop. 2.1, finishes the proof.
\end{proof}  
It is easily seen that $\TL$ has a subcategory $\TL^{ev}$
whose objects consist of even numbers of points, and with the same
morphisms between sets of even points as for $\TL$. The evaluation
$\TL^{ev}(\tau)$ is defined in complete analogy to $\TL(\tau)$.

\begin{corollary} If $\tau^2$ is a proper root of unity
of degree $\ell$ with $\ell$ odd, the negligible morphisms
form the unique non-zero proper ideal in 
$\TL^{ev}$.
\end{corollary}

\begin{proof}
If $\ell$ is odd, $p_{[\ell-1]}$ is a morphism in $\TL^{ev}$.
The proof of the last theorem goes through word for work (one only
needs to make sure that one stays within $\TL^{ev}$, which is easy
to check). 
\end{proof}

\vskip .5 in 
\v


\bigskip
\begin{thebibliography}{LLT}
\bibitem[F]{Freedman}, {\em A magnetic model with a possible Chern-Simons
phase}, preprint 2001, quant-ph/0110060.
\bibitem[GW]{pacific}  F. Goodman and H. Wenzl, The Temperley-Lieb
Algebra at Roots of Unity, Pac. J. Math., {\bf 161} (1993) 307-334.
\bibitem[J]{jones}  Jones, V. F. R. Index for subfactors. 
Invent. Math. 72 (1983), no. 1, 1--25.
\bibitem[W1]{sequence of projections}  H. Wenzl, On sequences of
projections,
Math. Rep. C.R. Acad. Sc. Canada {\bf 9} (1987) 5-9.
\bibitem[W2]{type A} H. Wenzl, Hecke algebras of type $A_n$ and
subfactors,
Invent. math. {\bf 92} (1988) 349-383.


\end{thebibliography}
\end{document}